\documentclass[12pt]{amsart}
\usepackage{amssymb,supertabular,graphicx}

\newcommand{\tbsp}{\rule{0pt}{18pt}}

\newtheorem{thm}{Theorem}[section]
\newtheorem{lemma}[thm]{Lemma}

\theoremstyle{definition}
\newtheorem{defn}[thm]{Definition}
\newtheorem{ex}[thm]{Example}

\theoremstyle{remark}
\newtheorem*{rmk}{Remark}
\newtheorem*{c1}{Case 1}
\newtheorem*{c2}{Case 2}
\newtheorem*{s1}{Step 1}
\newtheorem*{s2}{Step 2}
\newtheorem*{s3}{Step 3}

\newtheorem*{note}{Note}

\begin{document}

\title{Picard lattices of families of K$3$ surfaces}
\author[sarah-marie belcastro]{sarah-marie belcastro \\Department of Mathematics\\ University of
Northern Iowa\\ Cedar Falls, IA,  50614-0506.\\email:
\texttt{smbelcas@math.uni.edu} }
\maketitle
\thispagestyle{empty}

\begin{abstract}
Using toric geometry, lattice theory, and elliptic
surface techniques, we compute the Picard Lattice of certain K$3$
surfaces.  In particular, we examine the generic member of each of M.
Reid's list of 95 families of Gorenstein K$3$ surfaces which occur as
hypersurfaces in weighted projective 3-spaces.  As an application,
   we are able to determine whether the mirror family (in the sense of
mirror symmetry for K$3$ surfaces) for each one is also on Reid's list.

\end{abstract}

\pagestyle{empty}

\normalsize

\section{Introduction}\label{sec:beg}
 
	Given an integral convex polytope $\Delta$, we may associate to
$\Delta$ a toric variety $\mathbb P_\Delta$ and a family of hypersurfaces
$X_\Delta$ in $\mathbb P_\Delta$, referred to as \emph{toric
hypersurfaces}.  Generically, the members of the family are transversal
to the parent variety and have only singularities inherited from
$\mathbb P_\Delta$ (see
\cite{khov}).  The scope of this paper is limited to 2-dimensional toric
hypersurfaces whose parent varieties are weighted projective spaces.  

In 1979, M. Reid
classified and listed all families of K$3$ weighted projective
hypersurfaces with Gorenstein singularities (but did not publish this
work). Yonemura
\cite{yonemura} lists the weight-vectors for the associated weighted
projective spaces for each of the 95 families of surfaces. This list of
surfaces arises in many subfields of algebraic geometry, including
singularity theory and the birational geometry of three-folds.   The
present work was originally motivated by a question of Dolgachev about
K$3$ mirror symmetry on the ``Famous 95" (see Section \ref{sec:data} for
this application), but became focused on techniques for exhibiting and
analyzing elliptic fibrations on K3 surfaces (see Sections
\ref{sec:cdg} and
\ref{sec:comp} in particular).

\vspace{.2cm}

 We denote weighted projective space as
$\mathbb P(q_1,q_2,\dots,q_n)$, where the $q_i$ are the weights of
the projective variables.  In classifying the ``Famous 95", the
restrictions that the surfaces be K$3$ and Gorenstein determined which 
weight-vectors are possible.  Thus, each of the 95 families may be
defined by the weighted projective space in which it resides, as
follows.   Let $S$ be a generic member of one of the 95 families.

  The adjunction formula
in weighted projective space (see \cite{dolgwps}) induces the property
that
$\mathrm{deg}(S)  = \sum_i q_i =s$.  To construct the toric variety to
which each family is associated, we will use the standard toric geometry
correspondance between monomials $x^{\vec{p_i}}$ and integral lattice
points
${\vec{p_i}}$.  (Excellent references on toric varieties are 
\cite{fulton} and \cite{danilov}.)
We want the monomials of the family's equation to have degree $s$, or
equivalently, \newline we want our lattice points to lie in 
$\{\vec{x}
\in
\mathbb R^4
| \sum q_ix_i = s\}.$ To capture all such integer points (and only these
integer points), we consider the rational polytope $Q$ defined by the
convex hull 
$$Q = \mathrm{ conv}\{(s/q_1,\dots,0),\dots,(0,\dots,s/q_n)\}.$$
Generally, $Q$ is not integral.  Thus, we need only consider the
integral polytope $\Delta$ defined as the convex hull of all integral
points of $Q$, hereafter referred to as the
\emph{Newton polytope}. 

 Associated to $\Delta$ is a toric variety
$\mathbb P_\Delta$.  The equations describing a family of hypersurfaces
lying in
$\mathbb P_\Delta$ are obtained from the Laurent polynomials $\sum_j
\lambda_{j} \vec{x}^{\vec{p_j}}$ (see \cite{khov}).  Here $\vec{p_j}$ are
the integral points of $\Delta$, and as the $\lambda_j$ vary we obtain
different hypersurfaces in the family.

Computationally, it is useful to observe that  the condition
$\sum q_ix_i = s$ means that the rational polytope $Q$ lies in a
hyperplane of
$\mathbb R^4$. Using this condition, we may consider $Q$
as lying in $\mathbb R^3$ and do the convex hull computations
accordingly.  (Several algorithms exist for taking convex hulls; details
of this computation for the ``Famous 95" may be found in
\cite{sm}.)

Now recall that for any K$3$
surface $S$, $$H^2(S,\mathbb Z)
\cong (E_8)^2 \perp (U)^3 \text{ \cite[p.241]{bpv}}$$  where $U$ is the
hyperbolic plane $\begin{pmatrix} 0 & 1 \\ 1 & 0
\end{pmatrix}  \sim\begin{pmatrix} -2 & 1 \\ 1 & 0
\end{pmatrix}$.   We consider the
intersection form on $H^2(S,\mathbb Z)$, which gives us the bilinear form
on the lattice.    
 
  The group of linear equivalence classes of
Cartier divisors, denoted ${\rm Pic}(S)$,  injects into $H^2(S,\mathbb Z)$
for K$3$ surfaces
\cite[p.241]{bpv}.  The image of ${\rm Pic}(S)$ in $H^2(S,\mathbb Z)$
(equal to  $H^2(X,\mathbb Z) \cap
H^{1,1}(X,\mathbb C)$)  is the algebraic cycles in $H^2(S,\mathbb Z)$;
this has no torsion, so we may consider
${\rm Pic}(S)$ as a lattice, and call it the Picard lattice.
For ease of discussion, Section \ref{sec:latt} reviews the
lattice theory used later in the paper.

Computation of ${\rm Pic}(S)$ is nontrivial.  We are able to accomplish
it for the generic members of the ``Famous 95" families because of the
fortunate conjunction of two properties:  (1) they are hypersurfaces in
special weighted projective spaces, so we can determine the rank of
${\rm Pic}(S)$ (see Section \ref{sec:rhos}), and
(2) they are nondegenerate hypersurfaces in toric varieties, so we can
easily and explicitly desingularize them (see Section \ref{ssec:tc}). 
 We will
generally refer to the singular model as $\bar{S}$ and the smooth model
as
$S$.  Our method of computing
${\rm Pic}(S)$ is roughly as follows:  After determining $\rho(S)$, the
rank of
${\rm Pic}(S)$, we
apply toric techniques to desingularize
$\bar{S}$.  We use the nonsingular model with
$\rho(S)$ to exhibit an elliptic fibration. (All K$3$ surfaces with
$\rho(S) > 5$ are also elliptic surfaces, and most
singular K$3$ 2-dimensional toric hypersurfaces will fall into this
category or have singularities whose exceptional curves form elliptic
fibres.)  Section
\ref{sec:cdg} explains this process and contains theorems on the
existence of elliptic fibrations.   
Finally, we analyze the elliptic fibration to produce
${\rm Pic}(S)$. We introduce a number of techniques for analyzing a
fibration in Section \ref{sec:comp}.  

We end the paper with the application of these computations to mirror
symmetry, in Section \ref{sec:data}.  Table
\ref{ta:res1} lists ${\rm
Pic}(S)$ for each of the ``Famous 95" as well as the mirror lattice and
corresponding family, if any.  They are indexed by number as identified
in
\cite{yonemura}.


\section{Some Lattice Theory}\label{sec:latt}

This section summarizes general background on lattices.  Further
information  may be found in \cite{nikulin}.

\begin{defn}
 A lattice is a pair $(L,b)$ where
$L$ is a finite-rank free
$\mathbb Z$-module and $b$ is a $\mathbb Z$-valued nondegenerate symmetric
bilinear form on $L$.
\end{defn}
    
We will only
consider even lattices, i.e. those where $b(l,l)$ is even for all $l$. We
also denote $b(l,l)$ by $\langle l,l \rangle$.

\begin{defn}
 The \emph{discriminant} of a lattice is the
determinant of the matrix of the associated bilinear form.
\end{defn}

\subsection{Nikulin's results}\label{sec:nik}
 
Lattices are classified by their discriminant quadratic forms $q_L$
(often referred to simply as
$q$). The
discriminant form is defined on the discriminant group $G_L
= L^\star/L$, where
$L^\star =$ Hom$(L,\mathbb Z)$.     We define $q$ by $q_L(x) =
\langle x,x\rangle$ mod $2\mathbb Z$, where $x\in G_L$; note that $q$ is
the quadratic form associated to the bilinear form $b$ which defines the
lattice $L$.  We will define three classes of forms on $G_L$: 
$w^\epsilon_{p,k}, u_k, v_k$. An excellent explanation of these forms
can be found in
\cite{brieskorn}.

\begin{defn}
 The quadratic form
$w^\epsilon_{p,k}$ on
${\mathbb Z}_{p^k}$ is defined in cases, depending on the value of $p$.
\end{defn}
  
\begin{c1}
 For $p \neq 2$, $\epsilon \in
\{\pm 1\}$.  The form has generator value $q(1) = ap^{-k}$ (mod
$2{\mathbb Z}$).  For
$\epsilon = 1$ we choose
$a$ to be the smallest positive even number with a quadratic residue;
for
$\epsilon = -1$ we choose $a$ to be the
smallest positive even number without a quadratic residue.
\end{c1}
 
\begin{c2}
 For
$p = 2$, there are more possibilities.
For
$k = 1, \epsilon \in\{ \pm 1\}, w_{2,1}^\epsilon$ is defined as the
form with generator value
$q(1) = \epsilon /2$.  For $k \ge 1$ and $\epsilon \in \{
\pm 1, \pm 5\}, w_{2,k}^\epsilon$ is defined as the form with
generator value $q(1) =
\epsilon /2^k$. 
\end{c2}

\begin{defn}
 The forms $u_k,
v_k$ on
$\mathbb Z_{2^k}
\oplus
\mathbb Z_{2^k}$ are defined via their
matrices,which are $$u_k = 2^{-k}
\left(\begin{array} {cc} 0 & 1 \\ 1 & 0
\end{array}\right) \hspace{1cm} v_k = 2^{-k}\left(\begin{array} {cc}
2 & 1 \\ 1 & 2
\end{array}\right).$$  
\end{defn}

The following theorem is a combination of \cite[1.8.1,1.8.2]{nikulin}.

\begin{thm}

(i) Every nontrivial, nondegenerate
irreducible quadratic form on a finite abelian group is isomorphic to
one of $u_k, v_k, w_{p,k}^\epsilon$. \newline (ii) Every nondegenerate
quadratic form on a 
finite abelian group is isomorphic to an
orthogonal direct sum of $u_k, v_k, w_{p,k}^\epsilon$. \newline
(iii)  This representation of a quadratic form is not
unique.  
\end{thm}

\subsection{Example: $T_{p,q,r}$ and $M_{\vec{p},\vec{\iota},k}$
Lattices}\label{sec:oth}

\begin{figure}
\begin{center}
\includegraphics[scale=.8]{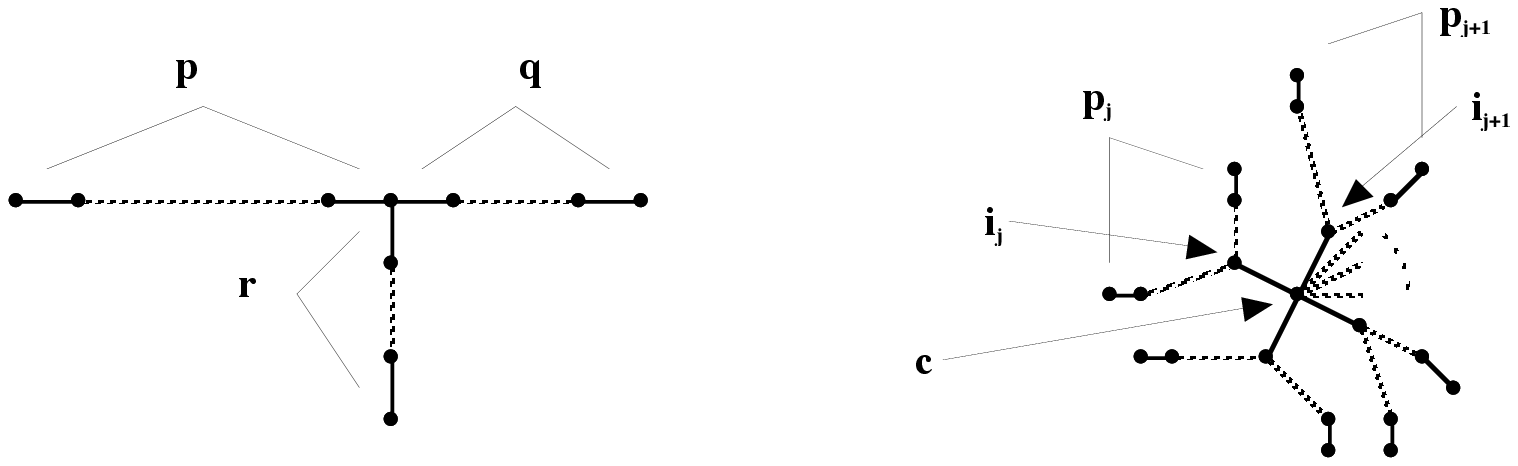}
\caption{$T_{p,q,r}$ and $M_{\vec{p},\vec{\iota},k}$ graphs}
\label{fig:tm}
\end{center}
\end{figure}

A $T_{p,q,r}$ lattice is so named because its graph forms the shape
of a
$\bf T$ (see Figure \ref{fig:tm}), where each vertex of the graph
represents a $(-2)$-curve.  In this notation,
$p$,$q$, and
$r$ are the lengths of the three legs, counting the central vertex each
time.  The discriminant of a $T_{p,q,r}$ is $pqr-pq-pr-qr$ and its rank
is
$p+q+r-2$.  Some of these lattices are isomorphic to Dynkin diagrams;
for example,
$T_{2,3,7}$ is
$E_8 \perp U$.  These lattices are well-known; they are carefully and
clearly studied in 
\cite{brieskorn}. 

A generalization of the $T_{p,q,r}$ lattices is 

\begin{defn}
Let $\vec{p} = (p_1,\dots,p_n)$ be an $n$-tuple of
positive integers, ordered from least to greatest. Let $\vec{i} =
(i_1,\dots,i_n)$ be an $n$-tuple of integers such that $i_j \leq
\lceil\frac{p_j}{2}\rceil$.  Let $k \geq -4$ be an even integer.
 Then $M_{\vec{p},\vec{\iota},k}$ is the lattice defined by the incidence
matrix of the following graph:  \newline Begin with a central vertex
$c$ with self-intersection $k$.  For each $j$, adjoin to this vertex a
Dynkin diagram of type $A_{p_j}$ by adding an edge between $c$ and
vertex $i_j$ of the Dynkin diagram.  Note that every vertex, except
possibly $c$, represents a $(-2)$-curve.
\end{defn}

 By taking the determinant of the associated matrix, Dolgachev has
calculated the discriminant of an
$M_{\vec{p},\vec{\iota},k}$; it is
$$\mathrm{disc}(M_{\vec{p},\vec{\iota},k}) = -(-1)^{\sum_{j=1}^n p_j}
(p_1+1)\ldots(p_n+1) (k+\sum_{j=1}^n\frac{i_j(p_j+1-i_j)}{p_j+1}).$$
When $k=-2, n = 3, i_j = 1$, we have a $T_{p,q,r}$ lattice.

Note that there is no algorithm for computing the associated
quadratic form to an $M_{\vec{p},\vec{\iota},k}$ (or even for a
$T_{p,q,r}$ -- see \cite{brieskorn}).

\section{Calculating $\rho(S)$}   \label{sec:rhos}

Henceforth in this paper, $S$ will denote the generic member of one of
the ``Famous 95" families.  We would like to determine
$\rho({S})$, the rank of the Picard lattice.  Let $\psi:
{S} \rightarrow \bar{S}$ be the resolution of singularities. We know that $\bar{S}$ has a
natural desingularization in terms of $A_{p_j}$ singularities (see
Section
\ref{ssec:tc}).  In this notation, the subscript refers to the
Milnor number of the singularity.

\begin{lemma}
\label{lem:gen}
 Let $\bar{S}$ be a generic surface in one
of the 95 families.  Then $\rho(\bar{S})=1$.
\end{lemma}

\begin{proof}[Proof of \ref{lem:gen}.]
 For each of the ``Famous 95"
families, the degree of the generic surface is  $s =
q_0 + q_1 + q_2 + q_3$ (again, by adjunction).
In \cite{cox}, the author examines the cup product 
$c :H^1(\bar{S},T_{\bar{S}})_0 \times H^{2,0}(\bar{S}) \rightarrow
H^{1,1}(\bar{S})$, where $H^1(\bar{S},T_{\bar{S}})_0$ represents the
variations in
$\bar{S}$ from 
 varying the coefficients in its equation ($T_{\bar{S}})$ is the
Zariski tangent space).  He concludes that if this map is surjective,
then for generic
$\bar{S}$,
$\mathrm{rk}({\rm Pic}(\bar{S})) = 1$.  We will now show that this map is
surjective when $\deg(\bar{S})=s$.  
Denote by $R=\bigoplus R_j$ the
Jacobian ring $\mathbb C [x_0,x_1,x_2,x_3]/\langle\vec{\partial
f}\rangle$, graded by degree.  In \cite[4.3]{dolgwps} the generators
of  $H_0^{2,0}(\bar{S})$ and $H_0^{1,1}(\bar{S})$ are computed.
$H_0^{1,1}(\bar{S})$ is generated by the monomials $m$ in $R$ such
that
$(1/(\deg (\bar{S})))(((\deg_{x_0} m)+1)q_0 + ((\deg_{x_1} m)+1)q_1 +
((\deg_{x_2} m)+1)q_2 + ((\deg_{x_3} m)+1)q_3) = 2$.  Because $\deg
\bar{S} = q_0 +q_1 + q_2 + q_3$, this expression picks out monomials
of total degree equal to $\deg (S)=s$.  Thus,  $H_0^{1,1}(\bar{S})  =
R_s$.  The same expression is used to find generators for
$H_0^{2,0}(\bar{S})$, except that the expression is set equal to 1
and thus only the monomial 1, of total degree zero, generates
$H_0^{2,0}(\bar{S})$ so that
$H_0^{2,0}(\bar{S}) = R_0$. Now, by definition (see
\cite[p.31]{bpv}), $H^1(\bar{S}, T_{\bar{S}})_0 = H^1(\bar{S},
\Omega_{\bar{S}})_0 = H_0^{1,1}(\bar{S})$ and so the map becomes 
$c: R_{s} \times H^{2,0}(\bar{S}) \rightarrow H^{1,1}(S)$.
$H^2(\bar{S},\mathbb C) = H^{2,0}(\bar{S}) \oplus H^{1,1}(\bar{S})
\oplus H^{0,2}(\bar{S})$ is torsion-free; therefore
$H^{2,0}(\bar{S}) = R_0$ and $H^{1,1}(\bar{S})  = R_s$, so that we may
interpret the map as
 $c: R_s \times R_0 \rightarrow R_s$. Because $c$ is
just the multiplication in $R$, it is an
isomorphism as $R_0$ is generated by 1.  Thus $c$
is surjective, and for generic
$\bar{S}$,
$\mathrm{rk}({\rm Pic}(\bar{S})) = 1$.\end{proof}

\begin{lemma}
\label{lem:ding}
 The contribution to $\mathrm{Pic}({S})$ from
desingularizing
$\bar{S}$ is
$\Sigma_j p_j$.
\end{lemma}

The proof of \ref{lem:ding} is well-known and thus left to the reader.

\begin{thm}
\label{thm:rk}
Let $S$ be the minimal nonsingular model of $\bar{S}$ and $p_j$ be the
types of the $A_{p_j}$ singularities of $\bar{S}$.
Then $\rho({S}) = 1 +
\Sigma_j p_j$.
\end{thm}
 
\begin{proof}[Proof of \ref{thm:rk}.] We combine \ref{lem:gen} and
\ref{lem:ding} to obtain the result.
\end{proof}

\begin{rmk}
It is interesting to note that there are few results
which compute the Picard number of the generic member of a family of
toric hypersurfaces.  If there were results for more general toric
hypersurfaces, one could apply most of the techniques listed later to a
wider class of objects.
\end{rmk}

\section{Finding Pic$(S)$ over $\mathbb Q$ }\label{sec:cdg}
 
\subsection{Desingularization}\label{ssec:tc}

To desingularize $\bar{S}$, we use the combinatorics of the Newton
polytope $\Delta$. Recall that singularities of $\bar{S}$ are inherited
from $\mathbb P_\Delta$.  We locate singularities of $\mathbb P_\Delta$
by examining the vertices and edges of $\Delta$ (see \cite{khov2} and
\cite{fulton}).  In fact, we may ignore vertices of $\Delta$ because
$\bar{S}$ does not intersect the corresponding torus orbits
(\cite{khov2}).  To check whether an edge $E_{ij}$ of $\Delta$ encodes a
singularity, we look at the dual cone generated by vectors $v_i, v_j$
normal to the incident faces $F_i, F_j$. (These vectors point to the
corresponding vertices of the polar dual polytope $\Delta^\circ$, which
are integral for K3 toric hypersurfaces (see \cite{bat}).)    The vectors
$v_i, v_j$ do not necessarily generate the lattice; we determine the
number of curves added
($n$, as in $A_n$) by the number of additional vectors needed to generate
the lattice (see \cite{fulton}, Ch. 2).

Additionally, each face of $\Delta$ corresponds to an irreducible
curve; its genus is determined by the number of integer points on the
face (see \cite{khov}, \cite{khov2}).  The structure of
$\Delta$ completely determines the intersections of the curves produced
by desingularization and those which arise from faces of $\Delta$.  For
example, any two faces $F_i, F_j$ of $\Delta$ intersect in an edge
$E_{ij}$ of
$\Delta$.  This corresponds geometrically to curves $C_{F_i}, C_{F_j}$;
their intersection number $C_{F_i}\cdot C_{F_j}$ is equal to the
number of lattice points on $E_{ij}$, minus one.  This follows from
\cite{fulton}, Ch. 5.  If there is a singularity at the intersection of
$C_{F_i}$ and $C_{F_j}$, then the intersection multiplicity gives the
number of copies of the chain(s) of
curves produced in resolving the singularities associated to
$E_{ij}$.   

We may use this information to create a graph, in which each curve is
represented by a vertex and where intersections between curves are
represented by edges.
The intersection  matrix of this graph (henceforth referred to as the
\emph{desingularization graph}) is a bilinear form, which is a lattice.
This lattice is of the same rank as Pic$(S)$, generated by curves on the
surface, and is certainly a sublattice of finite index of Pic$(S)$.
 
We now lay the groundwork for forming elliptic fibrations
from the desingularization graph.  Section
\ref{sec:comp} describes how to analyze these fibrations to obtain
Pic$(S)$.

\subsection{Forming Elliptic Fibrations}\label{sec:fibs}

\begin{defn}
An elliptic fibration is a regular map $\pi : S
\rightarrow B$ from a surface $S$ to some base curve $B$, such that the
general fibre
$\pi^{-1}(b)$ is an elliptic curve.
\end{defn}

  Because the
surfaces we study are K$3$, the base curve is isomorphic to $\mathbb P^1$. 
We only consider elliptic fibrations which have a finite number of
sections; under this condition, all sections are torsion sections and
are thus disjoint
\cite[Lemma, p.72]{miranda}.  There
may be finitely many reducible fibres, and the possible types were
listed by Kodaira \cite[V.7]{bpv}.

Thus, we'll be looking for subgraphs in the
desingularization graph which are isomorphic to graphs of fibres from
the Kodaira Classification.  Notice that because toric K$3$s only have
${A}_n$ singularities, we are limited to fibres which are
irreducible or which correspond to extended Dynkin diagrams.

\begin{thm}
\label{thm:fib}
 If $F$ is a configuration of curves on the
{\rm{K}}$3$ surface
$S$, identical to one of Kodaira's list of elliptic fibres, then
there exists an elliptic fibration
$\pi: S
\rightarrow
\mathbb P^1$ with general fibre linearly equivalent to $F$.
\end{thm}

\begin{proof}
 We see by using Riemann-Roch that $|F|$ has projective
dimension 1; using an argument similar to \cite[VIII.17.3]{bpv}, we can
show it has empty base locus.  Thus $\vert F\vert$
determines a map from $S$ to
$\mathbb P^1$.  Application of Bertini then shows that
the generic element of $\vert F\vert$ is smooth. 
\end{proof}

In practice, we want to partition each output graph into collections of
subgraphs corresponding to fibres, sections, and multisections.  In this
context, $U$ is the
intersection matrix of a section with an irreducible fibre. We begin by
finding a fibre, then mark all adjacent vertices as
sections/multisections, and partition the remaining subgraph into other
fibres and sections/multisections if possible.  However, sometimes we
will be left with a collection of vertices which do not form any
elliptic fibre.

\begin{lemma}
\label{thm:com}
 Suppose on the {\rm{K}}$3$ surface $S$ there exists
a fibration $\pi: S
\rightarrow
\mathbb P^1$  with exhibited fibre $F$.  Then, if a collection of
$-2$ curves disjoint from $F$ corresponds to a proper subset of an
extended Dynkin diagram, there exist curves sufficient to complete the
fibre. 
\end{lemma}

The proof of Lemma \ref{thm:com} is left to the reader.
 
\begin{note}
 In terms of partitioning output graphs, Lemma
\ref{thm:com} allows us to add vertices to subgraphs to form extended
Dynkin diagrams.
\end{note}

The effects of Lemma \ref{thm:com} can be seen in terms of linear
algebra as well.  For example, suppose the addition of only one curve is
necessary to complete a fibre.  The extended Dynkin diagrams
have positive semidefinite forms, so the determinants of their
intersection matrices are zero.  The intersection matrix associated to an
incomplete fibre has nonzero determinant.  When we complete the fibre by
adding one curve, the rank of the associated matrix (and thus
that of the fibration) does not increase because the determinant
becomes zero.

If there is more than one way to complete a fibre with the same
(minimal) number of curves, more information is necessary to decide
which fibre exists. The discriminant of the original intersection
matrix gives us a finite number of possibilities (see Section
\ref{sec:int}) for the discriminant of the fibration.  Usually comparing
the discriminant of the proposed fibration to
that of the original matrix is sufficient to determine which fibre we
have.

One of the most useful tools in analyzing possible elliptic fibrations
is the Shioda-Tate formula.

\begin{lemma}[\cite{shioda}]
 Let $f: S\rightarrow B$ be
an elliptic fibration (with section) of a smooth surface $S$, and let
$\rho$ be the rank of $\mathrm{Pic}(S)$.  Then
$$\rho = 2 + \sum_{F_i} ((\sharp\text{ components in }F_i)
- 1) + \mathrm{rk}(MW),$$ where  $\mathrm{rk}(MW)$ is the rank of the
Mordell-Weil group of sections, and $F_i$ ranges over all fibres  (note
that irreducible fibres will not contribute to the sum).  The number $2$
corresponds to  the contribution from a section and an irreducible
fibre.
\end{lemma}

  Most of the time we will find a fibration which shows that
$\mathrm{rk}(MW)=0$, i.e.
$MW$ is finite.    
When we have shown that
$\mathrm{rk}(MW)=0$, we have also shown that our fibration generates a
finite-index sublattice of
${\rm Pic}(S)$ because the rank is the same as that of ${\rm Pic}(S)$. 

Furthermore, we have 
\begin{lemma}
\label{lem:nam}   
$$| MW|^2 disc({\rm Pic}(S)) =
\prod_{F_i} disc(F_i).$$  {\rm This follows from}
\cite[Corollary1.7]{shioda}.
\end{lemma}

If we know that $MW$ is finite, we have additional information which
gives us an upper and lower bound on $| MW |$ = the number of
sections.  Our lower bound is the number of sections we have
exhibited in the fibration.  The upper bound is given by the gcd of
the orders of $G_{F_i}$, because the torsion subgroup of $MW$
(which in our case equals $MW$) embeds in the discriminant group of each
fibre
\cite[p.70]{miranda}.
 Often the bounds on  $| MW |$ and the use
of Lemma \ref{lem:nam} will be enough to tell us that
$MW$ is trivial (at which point we are done, as we have then shown that
the index in ${\rm Pic}(S) = 1$).

\begin{figure}[]
\normalsize

\begin{center}
{Number 26}
\end{center}
\begin{minipage}[t]{5.5cm}
\includegraphics[scale=.4]{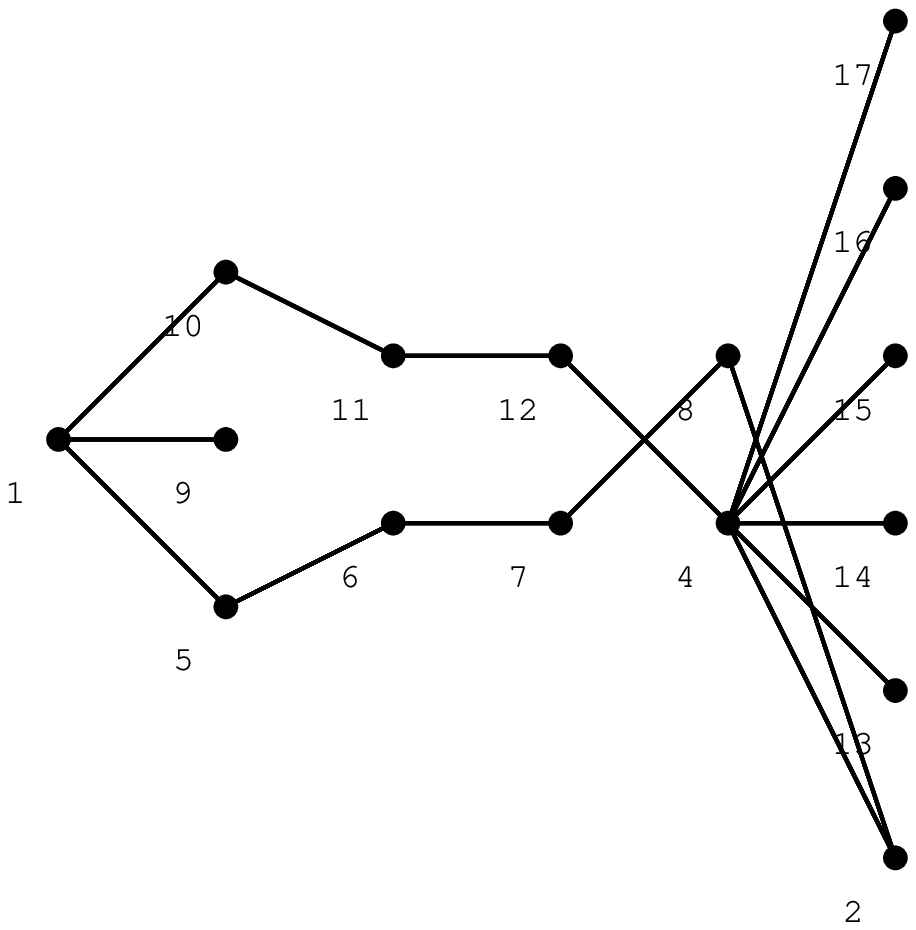}
\caption{\emph{Mathematica} output}
\label{fig:26i}
\end{minipage}
\hfill
\begin{minipage}[t]{5.5cm}
\includegraphics[scale=.4]{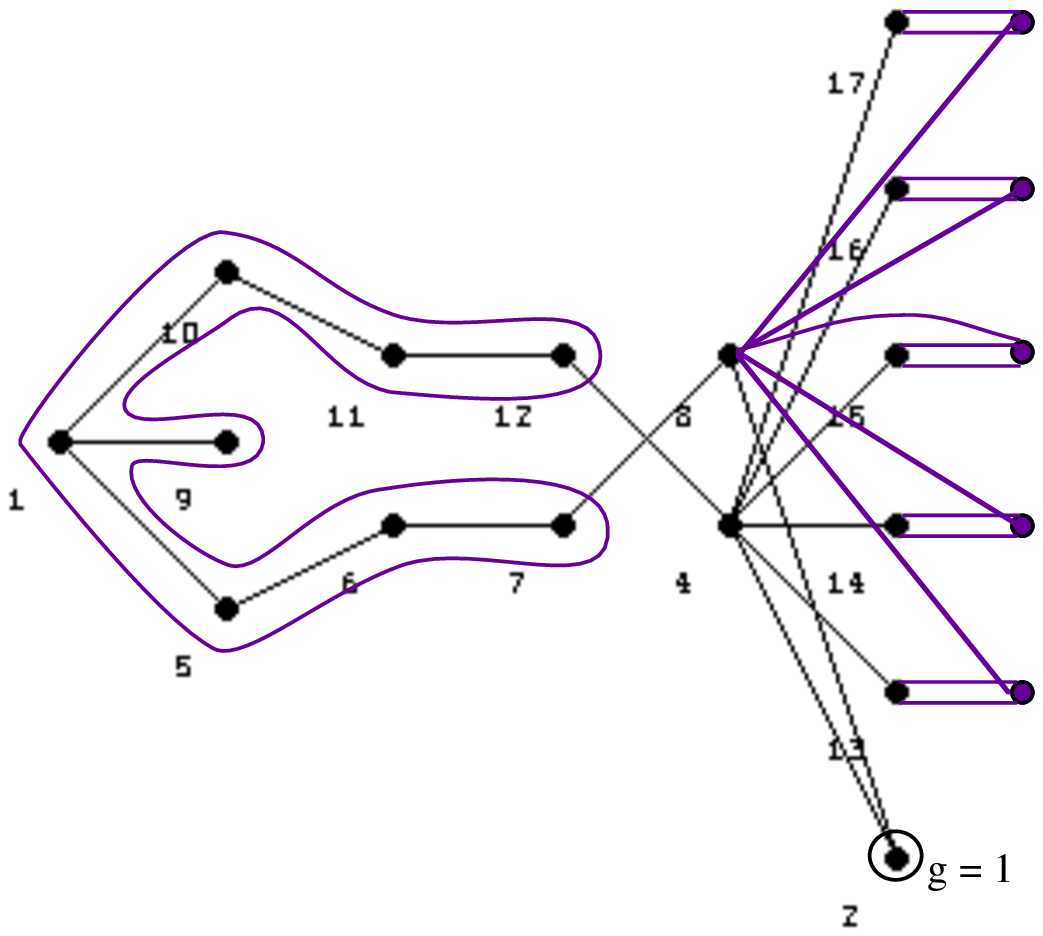}
\caption{Complete Fibration} \label{fig:26f}
\end{minipage}
\end{figure}

\begin{ex}\label{ex:1st}
 Figure \ref{fig:26i} shows the \emph{Mathematica}
output for surface family number 26.  Vertex 2 came from a face with one
interior point, and so is a curve of genus 1; by \ref{thm:fib}, it is
the generic fibre for an elliptic fibration.  We know from our
desingularization calculations (Theorem \ref{thm:rk}, Section
\ref{ssec:tc}) that $\rho=14$.
\end{ex}

Because curves 4 and 8 each intersect curve 2 exactly once, they must be
sections of the fibration.  Curves 1, 5-7, 9-12 are disjoint from curve
2, and form the reducible fibre
$\tilde{E}_7$ (see Figure \ref{fig:26f}).  This reducible fibre also
intersects each of curves 4 and 8 exactly once, as every section
intersects each fibre once.  Finally, using Theorem \ref{thm:com} we must
complete the remaining labeled curves 13-17 into a fibre or fibres.  We
know that each of these must be part of a different reducible fibre
because a section intersects each fibre with multiplicity one, and
section 4 intersects each of these curves with multiplicity one.  The
only fibre which conforms to these constraints is
$\tilde{A}_1$.  Furthermore, we can
only choose $\tilde{A}_1$ or we will violate the Shioda-Tate formula. 
Therefore, we add five vertices to the graph (see Figure
\ref{fig:26f}). 

The application of Shioda-Tate in this case gives us \newline
$14 = 2 + 7 + 5\cdot 1 +
\mathrm{rk}(MW) = 14 + \mathrm{rk}(MW)$ so $\mathrm{rk}(MW) = 0$.

Notice that because section 8 must intersect each fibre
once, and it does not naturally intersect any of curves 13-17, it must
intersect each of the added curves.

\section{Analyzing Fibrations to Compute ${\rm Pic}(S)$}\label{sec:comp}

\subsection{``Obvious" Elliptic Fibrations}
 \subsubsection{Elliptic Fibrations Begun with an Irreducible Fibre}

\begin{figure}[]
\normalsize

\begin{center}
{Number 65}
\end{center}
\begin{minipage}[t]{5cm}
\includegraphics[scale=.4]{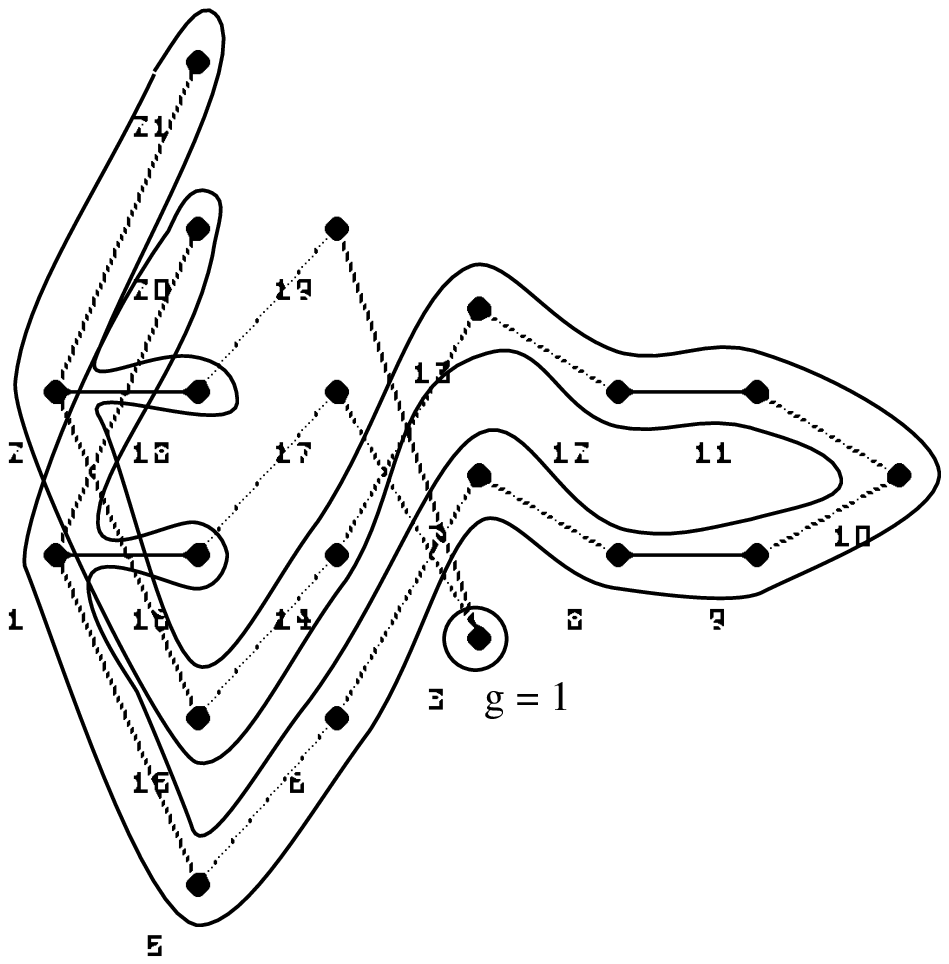}
\caption{First Fibration} \label{fig:651}
\end{minipage}
\hfill
\begin{minipage}[t]{5cm}
\includegraphics[scale=.4]{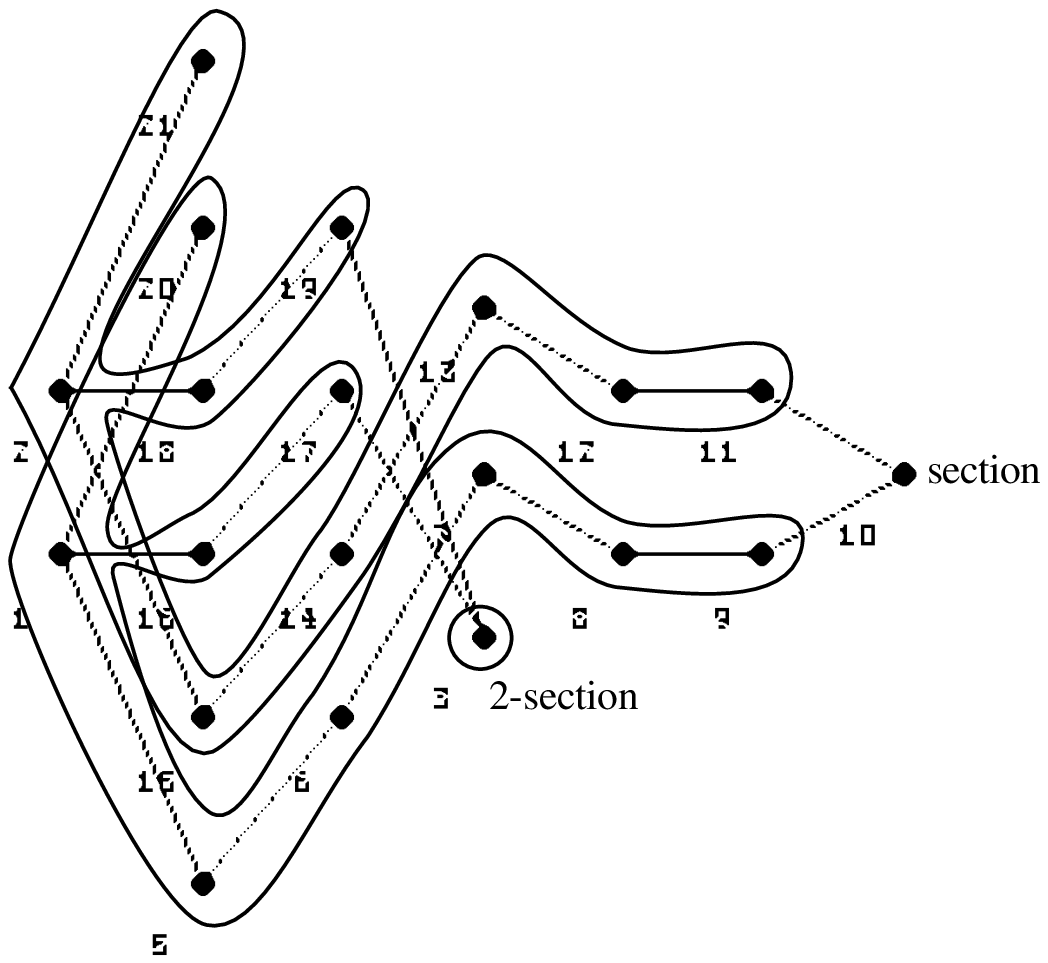}
\caption{Second Fibration} \label{fig:652}
\end{minipage}
\end{figure}

Figure \ref{fig:651} shows a fibration for number 65.  We are given via
our computer output that curve 3 has genus 1 and that $\rho = 18$. 
considering curve 3 as a fibre gives us curves 17 and 19 as
sections and the remaining curves as
$\tilde{D}_{16}$. 

Because 2 + 16 = 18, the Shioda-Tate formula
is satisfied with $\mathrm{rk}(MW) = 0$.    We have exhibited two
sections, and $G_{\tilde{D}_{16}} =
\mathbb Z_2
\oplus \mathbb Z_2$, so $2\leq | MW | \leq 4$.  As $| MW |^2
\in
\{4, 9, 16\}$, using
\ref{lem:nam} gives $| MW |^2 \mathrm{disc}({\rm
Pic}(S)) =4$.  Thus $| MW |^2 = 4$, so that $| MW | =2$ and
$\mathrm{disc}({\rm Pic}(S))=1$.  We notice immediately that therefore
${\rm Pic}(S) \neq D_{16} \perp U$ because that lattice has discriminant
4.  This indicates we might wish to look for another fibration, which
leads us to...

\subsubsection{Elliptic Fibrations with no Obvious Irreducible Fibre}

It happens frequently that we find more than one fibration.  For example,
in Figure
\ref{fig:652} we see another elliptic fibration for number 65, namely $E_8
\perp E_8 \perp U$.  

Application of Shioda-Tate shows that $\mathrm{rk}(MW)=0$, and formula
\ref{lem:nam} reads as
$1\cdot \mathrm{disc}({\rm Pic}(S)) =1$.  Thus 
${\rm Pic}(S) = E_8
\perp E_8 \perp U$.  Note that $D_{16}
\perp U$ has index 2 in $E_8
\perp E_8 \perp U$. In
\ref{sec:int}, we will show more generally how to determine ${\rm
Pic}(S)$ when its discriminant does not match that of the fibration.
 In cases where different fibrations produce seemingly different ${\rm
Pic}(S)$, we reconcile the discrepancy using the isomorphism relations
on the forms corresponding to the two lattices (see
\cite{brieskorn}, \cite{nikulin}).

\subsection{$T_{p,q,r}$ Fibrations} 
Sometimes we will not find any satisfactory fibrations using the
Dynkin diagrams, but will find a $T_{p,q,r}$ lattice (defined in
\ref{sec:oth}).

\begin{figure}[h]
\begin{center}
\includegraphics[scale=.5]{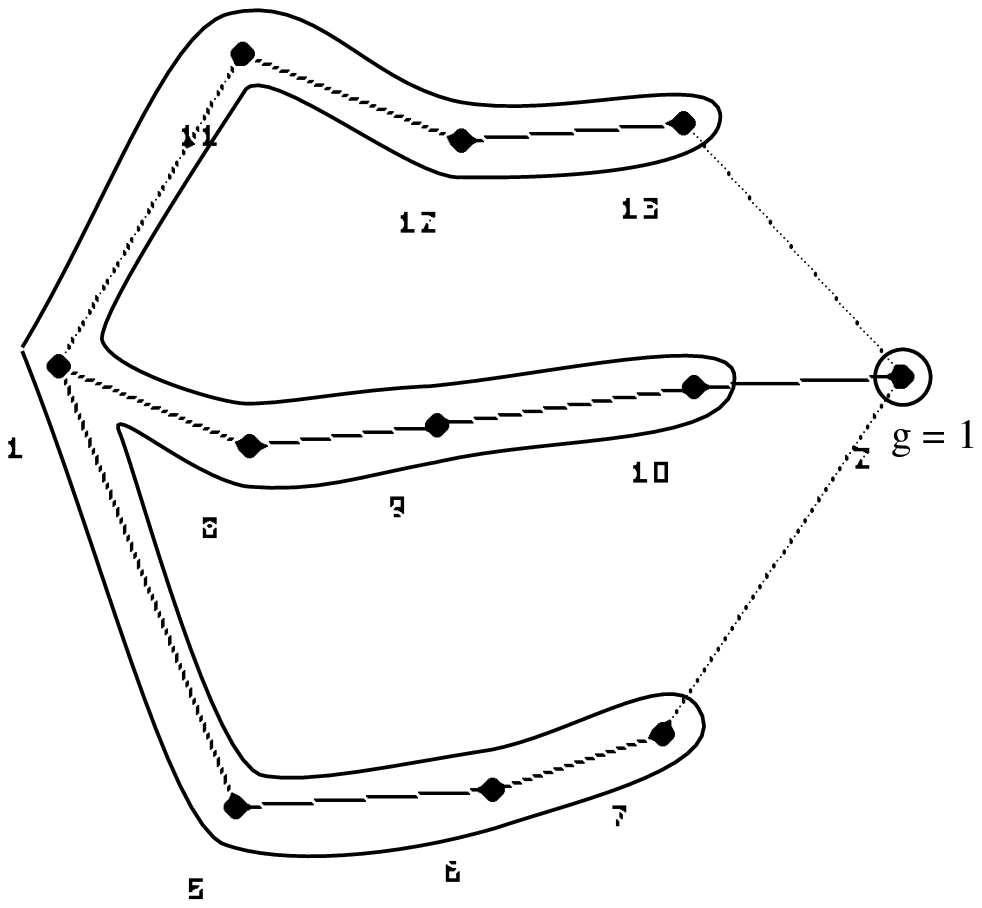}
\caption{Number 4} \label{fig:fib4}
\end{center}
\end{figure}

\begin{ex}
 Family number 4, which has $\rho = 10$, is shown in Figure
\ref{fig:fib4}.  The lattice $T_{4,4,4}$ exhibited in Figure
\ref{fig:fib4} has rank 10 as well.  Because number 4 is one of
Arnold's singularities, this is sufficient to compute ${\rm
Pic}(S)=T_{4,4,4}$, but for other cases we will need to use techniques
from sections
\ref{sec:int} and
\ref{sec:mat} to show that the $T_{p,q,r}$ lattice has index 1 in
Pic$(S)$.
\end{ex}

\subsection{Intermediate Lattice Calculations}\label{sec:int}

We will now resolve the uncertainties raised by the first fibration
(Figure \ref{fig:651}) for number 65.  More generally, we use
intermediate lattice calculations when we have exhibited more than one
section, especially when the number of exhibited sections divides
$\prod_{F_i} \mathrm{disc}(F_i)$.  First, recall that because we've
found a lattice of the correct rank, it must embed in
${\rm Pic}(S)$ with finite index. If $L$ is the  lattice
corresponding to our fibration,  then $L
\subseteq {\rm Pic}(S)
\subseteq L^\star$ (see Section \ref{sec:nik} for notation). 
Additionally, there is a 1-1 correspondence between the possible
``intermediate lattices"
$M$ and
$q$-isotropic subgroups $H$ of the discriminant group $G_L$
\cite[1.4.1(a)]{nikulin}.

\begin{figure}[h]
\begin{center}
\includegraphics[scale=.5]{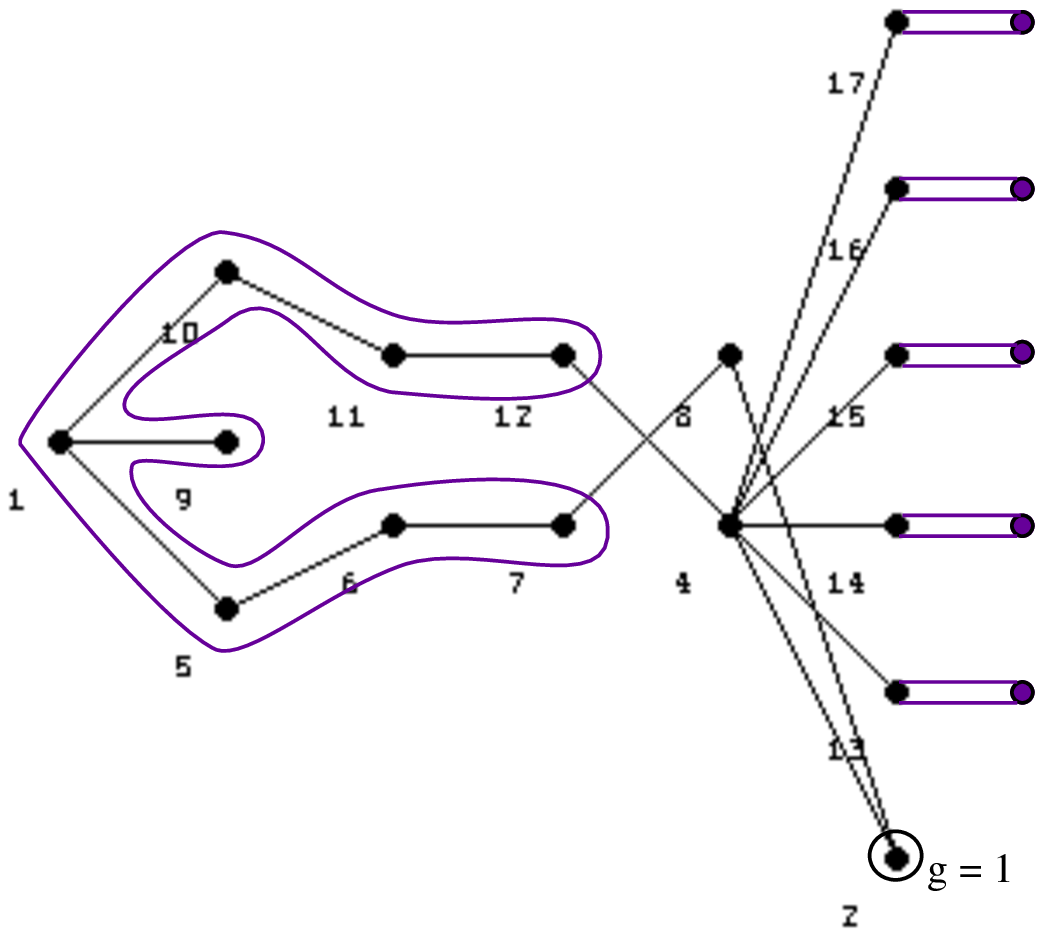}
\caption{Number 26} \label{fig:fib26}
\end{center}
\end{figure}

  In fact, there is
a constructive method for listing the different possibilities for
${\rm Pic}(S)$ via a formula of Nikulin:

\begin{thm}[\cite{nikulin}]
 For each $q_L$-isotropic
subgroup $H$ of $G_L$, $q_M = (q_L|H^\perp)/H$.
\end{thm}
  Let us interpret this
statement via an example.
\begin{ex}
 Number 26 (see Figure
\ref{fig:fib26}).
\end{ex}

\begin{s1}
 Find all $q_L$-isotropic subgroups of $G_L$.  A
subgroup is isotropic if every element of the subgroup has value
$q_L$-value 0 (mod 2$\mathbb Z$).
\end{s1}
   For number 26, we have
$\tilde{E_7} + 5\tilde{A_1}$.  Each of these degenerate fibres has
discriminant group $\mathbb Z_2$, so 
$G_L \cong(\mathbb Z_2)^6$.  The form corresponding to $E_7$ is $w_{2,1}^1$
and the value on the generator is 1/2.  The form corresponding to $A_1$
is $w_{2,1}^{-1}$ and the value on the generator is -1/2.  These
forms are independent, so in evaluating them on $(\mathbb Z_2)^6$ we can
just add the values on the components.  Immediate examples of
$q_L$-isotropic subgroups are those generated by (1,0,0,0,0,1),
(1,0,0,0,1,0), (1,0,0,1,0,0), (1,0,1,0,0,0), and (1,1,0,0,0,0). 
Fortunately, we 
 also know that 

\begin{thm}[\cite{nikulin}]
 Two intermediate lattices $M_1, M_2$
are isomorphic if and only if the corresponding $q_L$-isotropic subgroups
$H_1, H_2$ of $G_L$ are conjugate under an automorphism of $L$.
\end{thm}
 
\begin{rmk}
This corresponds, for example, to permuting several
copies of some Dynkin diagram, or equivalently to
permuting the relevant coordinates of several copies of some quadratic
form.
\end{rmk}
  So in our example, all five of the
$q_L$-isotropic subgroups listed above can be represented without
redundancy by (1,1,0,0,0,0).  In similar fashion, we have two other
distinct isotropic subgroups represented by (0,0,1,1,1,1) and
(1,1,1,1,1,1).

\vspace{.25cm}

\begin{s2}
 Determine $H^\perp$.
\end{s2} 
 To simplify the example, we
will only compute for $H$ = (0,0,1,1,1,1).
 Up
to permutation of the entries, we really only have 5 elements to deal
with: (0,0,0,0,0,0),
(0,0,0,0,0,1),
(0,0,0,0,1,1),
(0,0,0,1,1,1),
(0,0,1,1,1,1).
We will suppress the first 2 entries as they are always 0.  Now we
determine which of these are perpendicular to $H=\langle $(1,1,1,1)$\rangle$ with
respect to the quadratic form.  This holds true for an element $a$ when
$q(a) -q(a+$(1,1,1,1)$)=0$.

So respectively, for these 5 types of
elements, we have the data in Table \ref{ta:26h1p}.

\begin{table}[]
\normalsize

\begin{center}
\begin{tabular}{|l|l|l|}
\hline $a$ & $a+$(1,1,1,1) & $q(a) -q(a+$(1,1,1,1)$)$ \\
\hline (0,0,0,0) & (1,1,1,1) & 0 - 0 = 0 \\
\hline (0,0,0,1) & (1,1,1,0) & (-1/2) - (-3/2) = 1 \\
\hline (0,0,1,1) & (1,1,0,0) & (-1) - (-1) = 0 \\
\hline (0,1,1,1) & (1,0,0,0) & (-3/2) - (-1/2) = -1 \\
\hline (1,1,1,1) & (0,0,0,0) & 0 - 0 = 0 \\ \hline
\end{tabular}
\caption{Calculation of $H^\perp$ for Number 26} \label{ta:26h1p}
\end{center}
\end{table}

\begin{s3}
List all elements in $H^\perp$ and their values on
$q_L$.  Group them by conjugacy class in order to mod out by $H$. 
Using this list of values, determine the form of the intermediate
lattice corresponding to $H$.
\end{s3}
 Table \ref{ta:26h1c} lists all
elements of $H = \langle $(1,1,1,1)$\rangle$ by conjugacy class, and their
$q_L$-values. 

\begin{table}[]
\normalsize

\begin{center}
\begin{tabular}{|c|c|c|}
\hline elt & elt + (1,1,1,1) & $q_L$ value \\
\hline (0,0,0,0) & (1,1,1,1) & 0\\
\hline (0,0,1,1) & (1,1,0,0) & -1\\
\hline (0,1,0,1) & (1,0,1,0) & -1\\
\hline (0,1,1,0) & (1,0,0,1) & -1\\ \hline
\end{tabular} 
\caption{Values of the form $q_M$ corresponding to $H$ for Number 26}
\label{ta:26h1c}
\end{center}
\end{table}

This data corresponds to the form $v$.  We must also retain the original
form on the first two copies of $\mathbb Z_2$ ($w_{2,1}^{-1}
\perp w_{2,1}^1$) because they were not involved in the calculation;
they correspond to the zero-entries we suppressed above.

We have only computed $q_M$ for one of the
three distinct $q_L$-isotropic subgroups.  Number 26 is very
illustrative in that the other two intermediate lattices are $(w_{2,1}^{-1})^4$ and $u
\perp v$, neither of which is
isomorphic to $v \perp w_{2,1}^{-1}
\perp w_{2,1}^1$. This creates another question:  which one is
correct?   We need to look on the desingularization graph for other
fibrations which confirm that one of these choices is correct and that
the others are not possible. (It turns out that only $u\perp v$ is
possible; see
\cite[p.80]{sm}.)

Those hypersurface families for which we used intermediate lattice
calculations are 15, 16, 23, 26, 29, 30, 32, 34, 35, 46, 52, 54 - 56,
65, 68, 73 - 76, 80, 83, 84, 86, 92; details are in \cite{sm}.

\subsection{Methods for Fibrations Without Sections}\label{sec:muk}
Sometimes we'll only be able to find a fibration which has only
multisections, and no sections.    

\begin{thm}[\cite{ag2}]
\label{thm:jac}
Let  $\pi: S
\rightarrow B$ be an elliptic fibration.  Then there exists a fibration
$j:
{\bf J}_\pi(S) \rightarrow B$ with the following properties: \newline
$\bullet$ $j: {\bf J}_\pi(S) \rightarrow B$ has a section \newline
$\bullet$ each fibre of $\pi$ is isomorphic to some fibre of $j$
\newline
$\bullet$ if $\pi$ has a section, then ${\bf J}_\pi(S) \cong S$.
\end{thm}

We refer to $j:
{\bf J}_\pi(S) \rightarrow B$ in Theorem \ref{thm:jac} as
the \emph{Jacobian fibration}.  One way to construct ${\bf J}_\pi(S)$
from $S$ is to take the Jacobian variety of the generic
fibre $S_\eta$ and realize this as the generic fibre of some elliptic
surface which has a section.  By Theorem \ref{thm:jac}, the fibres of the
Jacobian fibration are the same as those of the non-Jacobian fibration;
likewise, $\rho(S) = \rho({\bf J}_\pi(S) )$.

\begin{figure}[]
\begin{center}
\includegraphics[scale=.5]{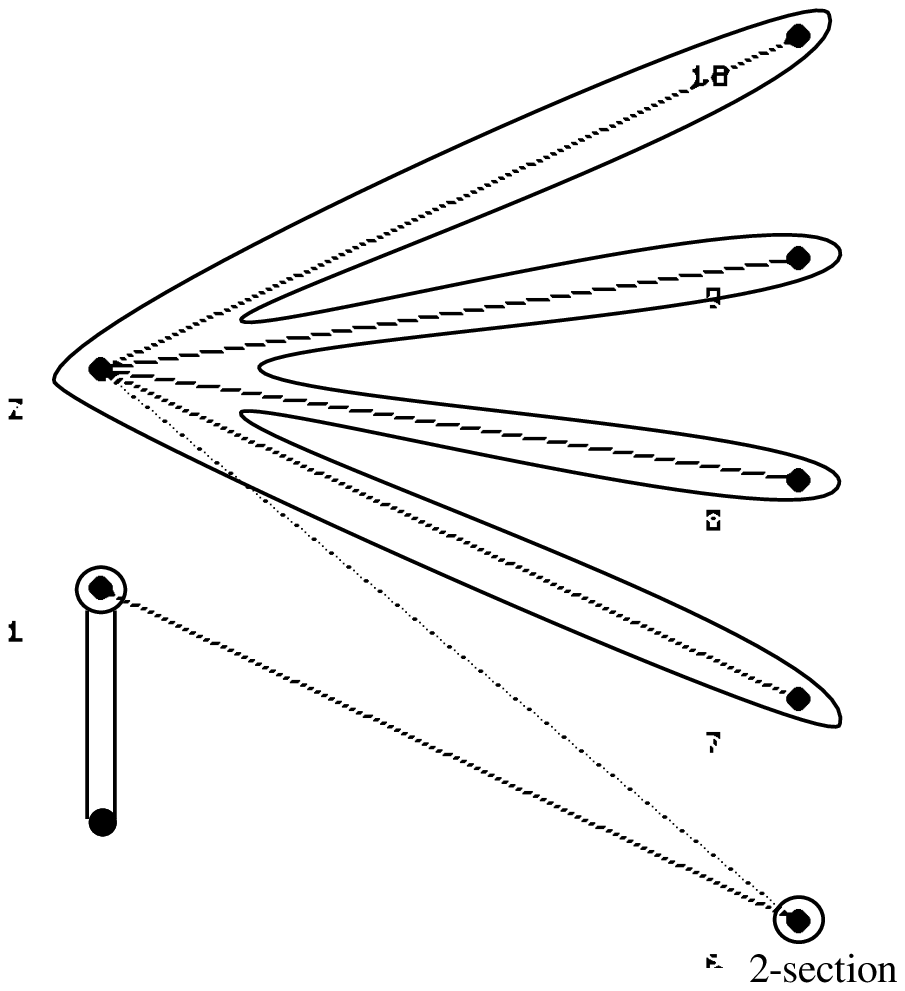}
\caption{Fibration for Number 19} \label{fig:fib19}
\end{center}
\end{figure}

\begin{ex}
In number 19 (Figure \ref{fig:fib19}) we immediately
see that curves 2, 7-10 form a $\tilde{D}_4$.  Thus curve 6 is a
2-section and we can add a curve, intersecting curve 1 twice, to form an
$\tilde{A}_1$.
\end{ex} 

We can show that this must be a non-Jacobian fibration using direct
computations on the intersection matrix.  
(In practice, we show by trial and error that if a section existed
for the fibration, then the rank of the intersection matrix would
exceed $\rho$.)
 For our example, number 19, there is
no fibration in which one can exhibit a section and satisfy Shioda-Tate
with $\mathrm{rk}(MW) = 0$.  
However, the Jacobian fibration associated to Figure \ref{fig:fib19}
satisfies Shioda-Tate with
$\mathrm{rk}(MW) = 0$;
$\rho = 7$ so we have \newline 7 = 2 + 4 + 1.

  In other words, a non-Jacobian fibration gives us information about
${\rm Pic}({\bf J}(S))$ and what we now need to know is the relationship
between
${\rm Pic}(S)$ and ${\rm Pic}({\bf
J}(S))$.    

We have shown (\cite{sm}) that there exists a
map
 $$\phi : {\rm Pic}(S)
\hookrightarrow {\rm Pic}({\bf J}(S))$$ of finite index
$n$, where $n$ is the index of multisections (the gcd of degrees of
multisections).  As ${\rm Pic}({\bf J}(S))$ is easy to compute, this map
along with application of lattice theory completely determine
${\rm Pic}(S)$. The current proof of
existence uses moduli spaces of sheaves; we are convinced that there is
a more direct proof and are working to complete one. This map was used to
calculate Pic$(S)$ in the nine cases for which we exhibited
non-Jacobian fibrations (numbers 2, 6, 18, 19, 31, 33, 53, 62, 69;
details in
\cite{sm}). In each case, we showed that  the existence of a
section increases  the Picard number.

\subsection{Calculating ${\rm Pic}(S)$ from the
Intersection Matrix}\label{sec:mat}

The intersection matrix for each desingularization graph represents a
quadratic form, which is calculable. This is particularly useful when
we cannot find a suitable fibration, or when we cannot show that
$MW$ is finite.
 Because the intersection matrix has rank equal to $\rho(S)$, it is of
finite index in
${\rm Pic}(S)$. Thus, if the discriminant of the matrix is square-free,
the index is 1 and we need merely determine the corresponding form.
(This was true for 
every case where 
$MW$ was of indeterminate rank.)  If the discriminant is not square-free,
then there are a finite number of possibilities for ${\rm Pic}(S)$ as
per Section \ref{sec:int}. 

\begin{note}
The matrices have dimension $(\rho +3) \times
(\rho +3)$, and rank $\rho$, so to determine the discriminant of such a
matrix we must find the minimum value of the determinants of the
$\rho\times\rho$ minors.
\end{note}
It should be mentioned that there is no general algorithm for
determining the quadratic form which corresponds to a matrix.  One
must decide based on the values of the form on its
generators and on the relations between these generators.  Those
surfaces for which we computed Pic$(S)$ from the matrix are 3, 7, 8, 17,
21, 36, 63, 66, 89, 94, 95; details are in \cite{sm}.

\section{Table of Picard Lattices and Application to Mirror
Symmetry}\label{sec:data}

Before presenting the available data on  the ``Famous
95," we briefly describe the original motivation for our work.

The 14 Arnold
singularities occur on Reid's ``Famous 95" list.  It is described
completely in \cite{dolgk3} that the ``strange duality" of Dolgachev and
Gabrielov for the 14 Arnold singularities is a K$3$ analogue to the more
recently studied Calabi-Yau threefold mirror symmetry.  The
question arises:  Can one find
 a larger class of surfaces which mirror each other?  The ``Famous 95"
are a natural set of K$3$ surfaces to investigate, as they arise so
often; Dolgachev remarked in \cite{dolgk3} that it would be most
interesting to see if all surfaces on the list mirror each other. 
(Surprisingly, they do not.\footnote{While Batyrev's mirror construction
does not always give the correct K$3$ mirror (see
\cite{dolgk3}), the spirit of his work suggests that  Kreuzer/Skarke's
 4319 toric K$3$ hypersurfaces
\cite{krsk} will mirror each other; the author is currently working to
confirm this.})  

\begin{defn}
 Two K$3$ surfaces form a mirror pair
$(S,\check{S})$ if 
$${\rm Pic}(S)^{\perp}_{H^2(S,\mathbb Z)}  = {\rm
Pic}(\check{S}) \perp U \mbox{ as lattices.}$$
\end{defn}

  We refer to
Pic$(\check{S})$  as the mirror lattice; to compute  Pic$(\check{S})$,
we use the fact that $L^{\perp}_{H^2(S,\mathbb Z)} = M$ if and only if
$q_L = -q_M$
\cite[1.6.2]{nikulin}.  

By computing ${\rm Pic}(S)$ for each of the ``Famous 95", we can see
which families mirror each other, and that the list is not
self-mirroring. The following table summarizes this data.   The ``No."
column indexes the surfaces as in \cite{yonemura}, as does the ``Mirror
Family" column, and the ``Weights" column gives the weights of the
projective variables in the corresponding weighted projective space.  As
for the notation in the
$\mathrm{Pic}(S)$ and ``Mirror Lattice" columns, the lattices
$A_n,D_n,E_n$ are the standard Dynkin lattices; $U$ is the hyperbolic
plane (see Section \ref{sec:beg}); the forms $u,v,$ and
$w^\epsilon_{p,k}$ are defined in Section
\ref{sec:nik}; and,
$T_{p,q,r}$ and
$M_{\vec{p},\vec{\iota},k}$ are defined in Section \ref{sec:oth}.   
Complete details of each calculation are in
\cite{sm} and available from the author.  


\begin{center}

\tablecaption{Pic$(S)$ for the 95 Families and Their
Mirrors}\label{ta:res1}
\tablefirsthead{\hline 
\multicolumn{1}{|c|}{\tbsp No.}
                       &
\multicolumn{1}{c}{Rank}
	&
\multicolumn{1}{c|}{Pic$(S)$}
                       &
\multicolumn{2}{c|} {Mirror Lattice and Family}
                      
& \multicolumn{1}{c|}{Weights} \\
\hline\tbsp  }
\tablehead{\hline \multicolumn{6}{|l|}{\small\sl
continued from previous page}\\
           \hline
\multicolumn{1}{|c|}{\tbsp No.}
                       &
\multicolumn{1}{c}{Rank}
	&
\multicolumn{1}{c|}{Pic$(S)$}
                       &
\multicolumn{2}{c|} {Mirror Lattice and Family}
                      
& \multicolumn{1}{c|}{Weights} \\

\hline\tbsp  }
\tabletail{\hline\multicolumn{6}{|r|}{\small\sl
continued on next
page}\\\hline}
\tablelasttail{\hline\hline}
\begin{supertabular}{| r | l |
p{3cm} | p{3cm} | l | l |}
1	&	$\rho = 1$ & $\langle 4\rangle$	&	\setlength{\rightskip}{0pt plus
1fil}$(E_8)^2 \perp \langle -4
\rangle
\perp U$ & 56, 73 &
$(1,1,1,1)$ \\ 
2	&	$\rho = 12$ & $E_6 \perp D_4 \perp U(3)$	&	$D_4
\perp A_2 \perp U(3)$	&	not on list	& $(2,3,3,4)$ \\ 
3 &	$\rho = 4$ & $M_{(1,1,1),(1,1,1),0}$	&	\setlength{\rightskip}{0pt
plus 1fil}$E_8\perp D_4\perp A_2
\perp U$	&	not on list	 & $(1,1,2,2)$ \\ 
4	&	$\rho = 10$ &
$T_{4,4,4}$ &	$T_{4,4,4}$	&	4	& $(1,3,4,4)$ \\ 
5	&	$\rho = 1$ &
$ \langle 2\rangle$	&
$(E_8)^2 \perp A_1 \perp U$	&	52, tetra.	& $(1,1,1,3)$ \\ 
6	&	$\rho
= 6$ & $D_4 \perp U(2)$	&	$D_8 \perp D_4 \perp U$	&	26, 34, 76	
 & $(1,2,2,5)$ \\ 
7	&	$\rho = 3$ & $M_{(1,1),(1,1),0}$	&
$E_8\perp D_7\perp U$	&	64	& $(1,1,2,4)$ \\ 
8	&	$\rho = 7$ &
$M_{(1,1,2,2),(1,1,1,1),-2}$	&	$q=w_{3,1}^1\perp w_{2,2}^{-1}$	&
not on list	& $(1,2,3,6)$ \\ 
9	&	$\rho = 10$ & $T_{2,5,5}$
	&	$T_{2,5,5}$	&	9,71	& $(1,4,5,10)$ \\ 
10	&
$\rho = 2$ & $U$	&	$(E_8)^2 \perp U$	&	65, 46, 80	 & $(1,1,4,6)$
\\ 
11	&	$\rho = 12$	&	$E_6 \perp D_4 \perp U$	&	$D_4 \perp A_2 \perp
U$	&	24	 & $(2,3,10,15)$ \\ 
12	&	$\rho = 6$	&	$D_4 \perp U$	&
$E_8 \perp D_4 \perp U$	&	27, 49	& $(1,2,9,6)$ \\ 
13	&	$\rho =
8$	&	$E_6 \perp U$	&	$E_8 \perp A_2 \perp U$	&	20, 59	&
$(1,3,8,12)$ \\ 
14	&	$\rho = 10$	&	$E_8 \perp U$	&	$E_8 \perp U$	&
14,28,45,51	& $(1,6,14,21)$ \\ 
15	&	$\rho = 14$	&	$E_6 \perp
(A_2)^3 \perp U$	&	$(A_2)^2 \perp U(3)$	&	not on list	 
& $(3,3,4,5)$ \\ 
16	&	$\rho = 16$	&	$E_8 \perp (A_2)^3 \perp U$	&
$A_2 \perp U(3)$	&	not on list	 & $(3,6,7,8)$ \\ 
17	&
$\rho = 14$	&	$T_{2,5,5}\perp A_4$	&\setlength{\rightskip}{0pt plus
1fil}	$A_4 \perp \left(\begin{array} {cc} 2 & 1 \\ 1 & -2
\end{array}\right)$ 	&	not on list	 & $(2,3,5,5)$ \\

$18$	&	$\rho = 8$	&	$M_{(1,2,2,2),(1,1,1,1),-2}$	&
$q = w_{3,2}^1 \perp w_{3,1}^1$	&	not on list	 & $(1,2,3,3)$
\\

$19$	&	$\rho = 7$	&
$\begin{array} {c}
M_{(1,1,1,1,2),} \\  _{(1,1,1,1,1),-2}
\end{array}$
&	$q =
v \perp w_{2,3}^1$	&	not on list	 & $(1,2,2,3)$ \\ 
20	&	$\rho =
12$	&	$E_8 \perp A_2 \perp U$	&	$E_6 \perp U$	&	13, 72	&
$(1,6,8,9)$ \\

21	&	$\rho = 2$	&	$\left(\begin{array} {cc}
2 & 1 \\ 1 & -2
\end{array}\right)$	&	$E_8 \perp T_{2,5,5}$	&	30, 86	& $(1,1,1,2)$ \\

22	&	$\rho = 10$	&	$E_6 \perp A_2 \perp U$	&	$E_6 \perp A_2 \perp U$
&	22	 & $(1,3,5,6)$ \\

23	&	$\rho = 11$	&	$D_5 \perp D_4 \perp U(2)$	&	$D_4 \perp A_3 \perp
U(2)$	&	not on list	 & $(2,2,3,5)$ \\

24	&	$\rho = 8$	&	$D_4 \perp A_2 \perp U$	&	$E_6 \perp D_4 \perp U$	&
11	& $(1,2,4,5)$ \\

25	&	$\rho = 4$	&	$A_2 \perp U$	&	$E_8 \perp E_6 \perp U$	&	43, 48, 88
& $(1,1,3,4)$ \\

26	&	$\rho = 14$	&	$D_8 \perp D_4 \perp U$ 	&	$D_4 \perp U(2)$	&	6	&
 $(2,4,5,9)$ \\

27	&	$\rho = 14$	&	$E_8 \perp D_4 \perp U$	&	$D_4 \perp U$	&	12	&
 $(2,3,8,11)$ \\

28	&	$\rho = 10$	&	$E_8 \perp U$	&	$E_8 \perp U$	&	14,28,45,51	
& $(1,3,7,10)$ \\

29	&	$\rho = 16$	&	$T_{2,5,5} \perp D_6$	&\setlength{\rightskip}{0pt
plus 1fil}	$q = w_{5,1}^{-1} \perp (w_{2,1}^{-1})^2$ 	&	not on list	 &
$(4,5,6,15)$ \\

30	&	$\rho = 18$	&	$E_8 \perp T_{2,5,5}$	&	$\left(\begin{array} {cc}
2 & 1 \\ 1 & -2
\end{array}\right)$
& 21 &	
 $(5,7,8,20)$ \\ 
$31$	&	$\rho = 15$	&	$E_6 \perp A_7 \perp
U$	&	$q = w_{2,3}^{-1} \perp w_{3,1}^1$	&	not on list	 &
$(3,4,5,12)$ \\ 
32	&	$\rho = 10$	&	$D_4 \perp D_4 \perp U(2)$	&	$D_4 \perp D_4 \perp
U(2)$	&	 32	 & $(2,2,3,7)$ \\ 
$33$	&	$\rho = 12$
&$\begin{array} {c}
M_{(1,1,1,1,2,2,3),} \\  _{(1,1,1,1,1,1,1),-4}
\end{array}$	 &
\setlength{\rightskip}{0pt plus 1fil}$q = w_{3,1}^1 \perp v \perp
w_{2,1}^1 \perp w_{2,1}^{-1}$	&	not on list	 & $(2,3,4,9)$ \\ 
34	&	$\rho = 14$	&	$D_8 \perp
D_4 \perp U$ 	&	$D_4 \perp U(2)$	&	6	 & $(2,6,7,15)$ \\

35	&	$\rho = 16$	&	$E_8 \perp A_6 \perp U$	&	$M_{(1,2),(1,1),0}$	&	66
&
 $(3,4,7,14)$ \\ 
36	&	$\rho = 13$	&	$T_{2,5,5}\perp A_3$
&	\setlength{\rightskip}{0pt plus 1fil}$D_5 \perp \left(\begin{array}
{cc} 2 & 1 \\ 1 & -2
\end{array}\right)$	&	not on list		 & $(2,3,5,10)$ \\

37	&	$\rho = 9$	&	$T_{3,4,4}$	&	$T_{2,5,6}$	&	58	 & $(1,3,4,8)$ \\

38	&	$\rho = 11$	&	$E_8 \perp A_1 \perp U$	&	$E_7 \perp U$	&	50, 82	&
 $(1,6,8,15)$ \\

39	&	$\rho = 9$	&	$E_6 \perp A_1 \perp U$	&	$E_7 \perp A_2 \perp U$	&
60	& $(1,3,5,9)$ \\

40	&	$\rho = 7$	&	$D_4 \perp A_1 \perp U$	&	$E_7 \perp D_4 \perp U$	&
81	 & $(1,2,4,7)$ \\

41	&	$\rho = 13$	&	$E_6 \perp D_5 \perp U$	&	$A_3 \perp A_2 \perp U$
&	not on list	 & $(2,3,7,12)$ \\

42	&	$\rho = 3$	&	$A_1 \perp U$	&	$E_8 \perp E_7 \perp U$	&	68, 83, 92
 & $(1,1,3,5)$ \\

43	&	$\rho = 16$	&	$E_8 \perp E_6 \perp U$	&	$A_2 \perp U$	&	25	
& $(3,4,11,18)$ \\

44	&	$\rho = 7$	&	$D_5 \perp U$	&	$E_8 \perp A_3 \perp U$	&	not on
list	& $(1,2,5,8)$ \\

45	&	$\rho = 10$	&	$E_8 \perp U$	&	$E_8 \perp U$	&	14,28,45,51	
& $(1,4,9,14)$ \\

46	&	$\rho = 18$	&	$E_8^2 \perp U$	&	$U$	&	10	 &
$(5,6,22,33)$ \\ 
47	&	$\rho = 15$	&	$E_7 \perp E_6 \perp U$	&	$A_2
\perp A_1 \perp U$	&	not on list	 & $(3,4,14,21)$ \\

48	&	$\rho = 16$	&	$E_8 \perp E_6 \perp U$	&	$A_2 \perp U$	&	25	
& $(3,5,16,24)$ \\

49	&	$\rho = 14$	&	$E_8 \perp D_4 \perp U$	&	$D_4 \perp U$	&	12	
& $(2,5,14,21)$ \\

50	&	$\rho = 9$	&	$E_7  \perp U$	&	$E_8 \perp A_1 \perp U$	&	38, 77	&
 $(1,4,10,15)$ \\

51	&	$\rho = 10$	&	$E_8 \perp U$	&	$E_8 \perp U$	&	14,28,45,51	
& $(1,5,12,18)$ \\

52	&	$\rho = 19$	&		\setlength{\rightskip}{0pt plus 1fil}$E_8 \perp
D_9 \perp U \cong (E_8)^2 \perp \langle -4\rangle \perp U$	&
$\langle 4\rangle, q=w_{2,2}^1$ & 5	  & $(7,8,9,12)$ \\

$53$	&	$\rho = 15$	&	$\begin{array} {c}
M_{(1,2,2,2,3,4),} \\  _{(1,1,1,1,1,1),-4}
\end{array}$		&	
\setlength{\rightskip}{0pt plus 1fil}$q = w_{3,2}^1 \perp w_{3,1}^1
\perp w_{2,1}^{-1}$	&	not on list	 & $(3,4,5,6)$ \\

54	&	$\rho = 16$	&	$E_8 \perp (A_2)^3  \perp U$	&	$A_2 \perp U(3)$	&
not on list	 & $(3,5,6,7)$ \\

55	&	$\rho = 15$	&	$D_9 \perp D_4 \perp U$	&	$w_{2,2}^1 \perp u$	&	not
on list	 & $(2,5,6,7)$ \\

56	&	$\rho = 19$	&	$E_8^2 \perp A_1 \perp U$	&	$\langle 2\rangle,
q=w_{2,1}^1$ & 1 &
$(5,6,8,11)$ \\

57	&	$\rho = 17$	&	\setlength{\rightskip}{0pt plus 1fil}$E_8
\perp D_5 \perp A_2 \perp U$	&	$w_{2,2}^5 \perp w_{3,1}^{-1}$	&	not on
list	& $(4,5,6,9)$ \\

58	&	$\rho = 11$	&	$T_{2,5,6}$	&	$T_{3,4,4}$	&	37	& $(1,4,5,6)$
\\

59	&	$\rho = 12$	&	$E_8 \perp A_2 \perp U$	&	$E_6 \perp U$	&	13, 72	&
 $(1,5,7,8)$ \\

60	&	$\rho = 11$	&	$E_7 \perp A_2 \perp U$	&	$E_6 \perp A_1 \perp U$
&	39	 & $(1,4,6,7)$ \\

61	&	$\rho = 18$	&	$E_8 \perp D_8 \perp U$	& $U(2)$	&	not on list	&
 $(4,6,7,11)$ \\

$62$	&	$\rho = 16$	&	$D_9 \perp D_5 \perp U$	&	$q = w_{2,2}^1
\perp w_{2,2}^5$	&	not on list	 & $(3,4,5,8)$ \\

63	&	$\rho = 8$	&	$M_{(1,1,2,3),(1,1,1,1),-2}$	&	$T_{2,5,5}
\perp (A_1)^2$	&	not on list	 & $(1,2,3,4)$ \\

64	&	$\rho = 17$	&	$E_8 \perp D_7 \perp U$	&	$M_{(1,1),(1,1),0}$	&	7	&
 $(3,4,7,10)$ \\

65	&	$\rho = 18$	&	$E_8^2 \perp U$	&	$U$	&	10	 &
$(3,5,11,14)$ \\

66	&	$\rho = 4$	&	$M_{(1,2),(1,1),0}$	&	$E_8 \perp A_6 \perp U$	&	35	&
$(1,1,2,3)$ \\

67	&	$\rho = 14$	&	\setlength{\rightskip}{0pt plus 1fil}	$E_6^2 \perp
U \cong E_8 \perp (A_2)^2 \perp U$	&	$(A_2)^2 \perp U$	&	not on list
& $(2,3,7,9)$ \\

68	&	$\rho = 17$	&	$E_8 \perp E_7 \perp U$	&	$A_1 \perp U$	&	42	& 
 $(3,4,10,13)$ \\

$69$	&	$\rho = 13$	&	$D_4 \perp A_7 \perp U$	&	$q =
w_{2,3}^{-1}
\perp v$	&	not on list	 & $(2,3,4,7)$ \\

70	&	$\rho = 14$	&	\setlength{\rightskip}{0pt plus 1fil}$E_8 \perp
A_2 \perp (A_1)^2 \perp U$	&\setlength{\rightskip}{0pt plus 1fil}	$q=
w_{3,1}^{-1} \perp (w_{2,1}^1)^2$	& not on list	  & $(2,3,5,8)$ \\

71	&	$\rho = 10$	&	$T_{2,5,5}$	&	$T_{2,5,5}$	&	9, 71	 &
$(1,3,4,7)$ \\

72	&	$\rho = 8$	&	$E_6 \perp U$	&	$E_8 \perp A_2 \perp U$	&	20, 59	&
 $(1,2,5,7)$ \\

73	&	$\rho = 19$	&	$E_8^2 \perp A_1 \perp U$	&	$\langle 2\rangle,
q=w_{2,1}^1$ & 1 &	
 $(7,8,10,25)$ \\

74	&	$\rho = 17$	&	$M_{(3,3,4,6),(1,1,1,3),-4}$	&	 $q =
w_{2,3}^{-5}$ & not on list	 & $(4,5,7,16)$ \\

75	&	$\rho = 13$	&	$E_7 \perp (A_1)^4 \perp U$	&	$(A_1)^5 \perp U$ 	&
not on list	 & $(2,4,5,11)$ \\

76	&	$\rho = 14$	&	$D_8 \perp D_4 \perp U$	&	$D_4 \perp U(2)$	&	6	&
 $(2,5,6,13)$ \\

77	&	$\rho = 11$	&	$E_8 \perp A_1 \perp U$	&	$E_7 \perp U$	&	50, 82	&
 $(1,5,7,13)$ \\

78	&	$\rho = 10$	&	$E_7 \perp A_1 \perp U$	&	$E_7 \perp A_1 \perp U$
&	78  & $(1,4,6,11)$ \\

79	&	$\rho = 15$	&	$E_8 \perp D_5 \perp U$	&	$A_3 \perp U$	&	not on
list	 & $(2,5,9,16)$ \\

80	&	$\rho = 18$	&	$E_8^2 \perp U$	&	$U$	&	10	&
$(4,5,13,22)$ \\

81	&	$\rho = 13$	&	$E_8 \perp (A_1)^3 \perp U$	&	$D_4 \perp A_1 \perp
U$ &	40	 & $(2,3,8,13)$ \\

82	&	$\rho = 9$	&	$E_7 \perp U$	&	$E_8 \perp A_1 \perp U$	&	38, 77	&
 $(1,3,7,11)$ \\

83	&	$\rho = 17$	&	$E_8 \perp E_7 \perp U$	&	$A_1 \perp U$	&	42	&	
$(4,5,18,27)$ \\
84	&	$\rho = 18$	&	$E_8 \perp A_8 \perp U$ 	& $q =
w_{3,2}^{-1}$ & not on list	& $(5,6,7,9)$ \\
85	&	$\rho = 13$	&\setlength{\rightskip}{0pt plus 1fil}	$D_4 \perp A_6
\perp A_1 \perp U$	&	$D_4 \perp
\langle 14\rangle \perp U$ 	&	not on list	 & $(2,3,4,5)$ \\
86	&	$\rho = 18$	&	$E_8 \perp T_{2,5,5}$	&	$\left(\begin{array} {cc}
2 & 1 \\ 1 & -2
\end{array}\right)$	&	21	 & $(4,5,7,9)$ \\
87	&	$\rho = 10$	&	$T_{3,4,5}$	&	$T_{3,4,5}$	&	87	 & $(1,3,4,5)$
\\
88	&	$\rho = 16$	&	$E_8 \perp E_6 \perp U$	&	$A_2 \perp U$	&	not on
list	 & $(2,5,9,11)$ \\
89	&	$\rho = 8$	&	$M_{(1,2,4),(1,1,2),-2}$	&	$A_{10}\perp U$	&	not on
list	&
 $(1,2,3,5)$ \\

90	&	$\rho = 17$	&\setlength{\rightskip}{0pt plus 1fil}	$E_8 \perp D_6
\perp A_1 \perp U$	&	$A_1 \perp U(2)$	&	not on list	 & $(4,6,7,17)$ \\

91	&	$\rho = 18$	&	$E_8 \perp E_7 \perp A_1 \perp U$	&	$q=w_{2,1}^1
\perp w_{2,1}^{-1}$	&	not on list	 & $(5,6,8,19)$ \\
92	&	$\rho = 17$	&	$E_8 \perp E_7 \perp U$	&	$A_1 \perp U$	&	42	&	
$(3,5,11,19)$ \\
93	&	$\rho = 16$	&	$E_8 \perp D_6 \perp U$	&	$(A_1)^2 \perp U$	&	not
on list	  & $(3,4,10,17)$ \\
94	&	$\rho = 16$	&	$M_{(2,3,4,6),(1,1,2,2),-4}$	&	$q=w_{19,1}^1$	&
not on list &	 $(3,4,5,7)$ \\
95	&	$\rho = 14$	&	$M_{(1,2,4,6),(1,1,2,3),-4}$	&	$q=w_{17,1}^{-1}$	&
not on list & $(2,3,5,7)$ \\
 \end{supertabular}
\end{center}

\vspace{.5cm}

\begin{center}\textsc{Acknowledgements} \end{center}
 The author is grateful to I. Dolgachev for
supervising this work, which was part of her Ph.D. dissertation at the
University of Michigan.  She would also like to thank   W. Cherry for
general consultations,  J.H. Keum for pointing out an omission in a
previous version, and  R. Miranda and D. Cox for providing assistance
via email.

\end{document}